\documentstyle[11pt]{article}
\input amssym.def
\input amssym.tex

\font\bbb=msbm10 

\newtheorem{thm}{Theorem}[section]
\newtheorem{defn}[thm]{Definition}

\newtheorem{prop}[thm]{Proposition}
\newtheorem{cor}[thm]{Corollary}

\newtheorem{rema}[thm]{Remark}

\newcommand{\halmos}{\rule{1ex}{1.4ex}}
\newcommand{\pfbox}{\hspace*{\fill}\mbox{$\halmos$}}
\newcommand{\Z}[0]{\mbox{\bbb Z}_+}
\newcommand{\Dz}[0]{\frac{\partial}{\partial \theta} + \theta
\frac{\partial}{\partial z}} 
\newcommand{\Lz}[0]{z^{j + 1} \frac{\partial}{\partial z} + (\frac{j +
1}{2})\theta z^j \frac{\partial}{\partial \theta}}
\newcommand{\Gz}[0]{z^j \left( \frac{\partial}{\partial \theta} -
\theta \frac{\partial}{\partial z}\right)}
\newcommand{\Loz}[0]{z \frac{\partial}{\partial z} +
(\frac{1}{2})\theta \frac{\partial}{\partial \theta}}
\newcommand{\Lminusz}[0]{z^{- j + 1} \frac{\partial}{\partial z} +
(\frac{- j + 1}{2})\theta z^{-j} \frac{\partial}{\partial \theta}}
\newcommand{\Gminusz}[0]{z^{- j + 1}\left( \frac{\partial}{\partial
\theta} - \theta \frac{\partial}{\partial z}\right)}

\begin{document}

\renewcommand{\theequation}{\thesection.\arabic{equation}}

\begin{center}
\begin{LARGE}
{\bf A supergeometric interpretation of vertex operator superalgebras}
\\ 
\end{LARGE}
\medskip
Katrina Barron \\
Rutgers University \\
\end{center}
\bigskip

\section{Introduction}

Conformal field theory (or more specifically, string theory) and
related theories (cf. \cite{BPZ}, \cite{FS}, \cite{V}, and \cite{S})
are the most promising attempts at developing a physical theory that
combines all fundamental interactions of particles, including gravity.
The geometry of this theory extends the use of Feynman diagrams,
describing the interactions of point particles whose propagation in
time sweeps out a line in space-time, to one-dimensional ``particles''
(strings) whose propagation in time sweeps out a two-dimensional
surface.  For genus zero holomorphic conformal field theory,
algebraically, these interactions can be described by products of
vertex operators or more precisely, by means of vertex operator
algebras (cf. \cite{Bo} and \cite{FLM}).  However, until 1990 a
rigorous mathematical interpretation of the geometry and algebra
involved in the ``sewing'' together of different particle
interactions, incorporating the analysis of general analytic
coordinates, had not been realized.  In \cite{H1} and \cite{H2},
motivated by the geometric notions arising in conformal field theory,
Huang gives a precise geometric interpretation of the notion of vertex
operator algebra by considering the geometric structure consisting of
the moduli space of genus zero Riemann surfaces with punctures and
local coordinates vanishing at the punctures, modulo conformal
equivalence, together with the operation of sewing two such surfaces,
defined by cutting discs around one puncture from each sphere and
appropriately identifying the boundaries.  Important aspects of this
geometric structure are the concrete realization of the moduli space
in terms of exponentials of a representation of the Virasoro algebra
and a precise analysis of sewing using these resulting exponentials.
Using this geometric structure, Huang then introduces the notion of
geometric vertex operator algebra with central charge $c \in \Bbb{C}$,
and proves that the category of geometric vertex operator algebras is
isomorphic to the category of vertex operator algebras.

In \cite{F}, Friedan describes the extension of the physical model of
conformal field theory to that of superconformal field 
theory and the notion of a superstring whose propagation in time
sweeps out a supersurface.  Whereas conformal field theory attempts to
describe the interactions of bosons, superconformal field theory
attempts to describe the interactions of boson-fermion pairs.  This,
in particular, requires an operator $D$ such that $D^2 =
\frac{\partial}{\partial z}$.  Such an operator arises naturally in
supergeometry.  In \cite{BMS}, Beilinson, Manin and Schechtman study
some aspects of superconformal symmetry, i.e., the Neveu-Schwarz
algebra, {}from the viewpoint of algebraic geometry.  In this work, we
will take a differential geometric approach, extending Huang's
geometric interpretation of vertex operator algebras to a
supergeometric interpretation of vertex operator superalgebras.   

Within the framework of supergeometry (cf. \cite{D}, \cite{R} and
\cite{CR}) and motivated by superconformal field theory, we define the
moduli space of super-Riemann surfaces with genus zero ``body'',
punctures, and local superconformal coordinates vanishing at the
punctures, modulo superconformal equivalence.  We announce the result
that any local superconformal coordinates can be expressed in terms of
exponentials of certain superderivations, and that these
superderivations give a representation of the Neveu-Schwarz algebra
with zero central charge.  We define a sewing operation on this moduli
space and give an interpretation of sewing in terms of these
exponentials of representatives of Neveu-Schwarz algebra elements.  We
then introduce the notion of {\it supergeometric vertex operator
superalgebra with central charge $c \in \Bbb{C}$}.  The purpose of
this paper is to announce the result that the category of
supergeometric vertex operator superalgebras with central charge $c
\in \Bbb{C}$ is isomorphic to the category of (superalgebraic) vertex
operator superalgebras with central charge $c \in \Bbb{C}$,
appropriately defined.

Recall that in a vertex operator algebra, the Virasoro element $L(-1)$
plays the role of the differential operator $\frac{\partial}{\partial
z}$.  Thus in considering what should be the corresponding
superalgebraic setting for superconformal field theory, we naturally
want to consider a super-extension of the Virasoro algebra, namely the
Neveu-Schwarz algebra \cite{NS} in which the element $G(-\frac{1}{2})$
has the supercommutator $\frac{1}{2} [G(-\frac{1}{2}),
G(-\frac{1}{2})] = L(-1)$.  In this work, we will assume that a vertex
operator superalgebra (cf. \cite{T}, \cite{G}, \cite{FFR}, \cite{DL},
and \cite{KW}) includes the Neveu-Schwarz algebra by definition (cf. 
\cite{KW}), but in addition, we extend the notion of vertex 
operator superalgebra to be over a Grassmann algebra instead of just
{\bbb C} and to include ``odd'' formal variables instead of just even
formal variables.  This notion of {\it vertex operator superalgebra
over a Grassmann algebra and with odd formal variables and central
charge $c \in \Bbb{C}$} is in fact equivalent to the notion of vertex
operator superalgebra over a Grassmann algebra without odd formal
variables.  However, in a vertex operator superalgebra with odd
variables, the fact that $G(-\frac{1}{2})$ plays the role of the
operator $D$ (mentioned above in reference to the supergeometry) is
made explicit and the correspondence with the supergeometry is more
natural.

The main result we announce is that the category of supergeometric
vertex operator superalgebras over a Grassmann algebra $\Lambda_*$
with central charge $c \in \Bbb{C}$ and the category of vertex
operator superalgebras over $\Lambda_*$ with central charge $c \in
\Bbb{C}$ and with (or without) odd formal variables are isomorphic.
Details of the proof of this result can be found in \cite{B}.

\section{Superconformal superfunctions and the Neveu-Schwarz algebra}

In this section, we follow many of the conventions developed in the
theory of superfunctions (cf. \cite{D}, \cite{R}).  Let
$\Lambda_\infty$ be the Grassmann algebra over {\bbb C} on an infinite
number of generators $\zeta_1, \zeta_2,...$, and let $I_\infty = \{(i)
= (i_1, i_2, \ldots, i_{2n}): i_1 < i_2 < \cdots < i_{2n}, i_l \in \Z,
n \in \Bbb{N}\}$, $J_\infty = \{(j) = (j_1, j_2, \ldots, j_{2n+1}):
j_1 < j_2 < \cdots < j_{2n + 1}, j_l \in \Z, n \in \Bbb{N} \}$, and
$K_\infty = I_\infty \cup J_\infty$. As a vector space
$\Lambda_\infty$ has a natural $\Bbb{Z}_2$-grading given by
$\Lambda_\infty = \Lambda_\infty^0 \oplus \Lambda_\infty^1$ where 
$\Lambda_\infty^0 = \{a \in \Lambda_\infty : a = \sum_{(i) \in
I_\infty} c_{(i)} \zeta_{i_1} \zeta_{i_2} ... \zeta_{i_{2n}} , \;
c_{(i)} \in \Bbb{C} \}$ is the {\it even} subspace and
$\Lambda_\infty^1 = \{a \in \Lambda_\infty : a = \sum_{(j) \in
J_\infty} c_{(j)} \zeta_{j_1} \zeta_{j_2} ... \zeta_{j_{2n + 1}} , \;
c_{(j)} \in \Bbb{C} \}$ is the {\it odd} subspace. (Note that
$\zeta_{(\emptyset)} = 1$.)  We can also decompose $\Lambda_\infty$
into {\it body} $(\Lambda_\infty)_B = \{c_{(\emptyset)} \in \Bbb{C}
\}$ and {\it soul} $(\Lambda_\infty)_S = \{a \in \Lambda_\infty : a =
\sum_{ {(k) \in K_\infty}\atop{k \neq (\emptyset)}} c_{(k)}
\zeta_{k_1} \zeta_{k_2} ... \zeta_{k_n} , \; c_{(k)} \in \Bbb{C} \}$
subspaces such that $\Lambda_\infty = (\Lambda_\infty)_B \oplus
(\Lambda_\infty)_S$.  For $a \in \Lambda_\infty$, we write $a = a_B +
a_S$ for its body and soul decomposition.

Let $z_B$ be a complex variable and $h(z_B)$ analytic.  For $z$ a
variable in $\Lambda_\infty^0$, we define $h(z) = \sum_{n \in \Bbb{N}}
\frac{z_S^n}{n!} h^{(n)}(z_B)$.  Note that if $h(z_B)$ is convergent
in an open neighborhood $N \subseteq \Bbb{C}$ of $z_B$ then $h(z)$ is
well-defined (i.e., convergent) in the open neighborhood $\{z = z_B +
z_S \in \Lambda_\infty^0 : z_B \in N \} \subseteq \Lambda_\infty$.
Let $f(z) = \sum_{(k) \in K_\infty} f_{(k)} (z) \zeta_{k_1}
\zeta_{k_2} \cdots \zeta_{k_n}$ where each $f_{(k)} (z_B)$ is
analytic.  We say that $f$ is a {\it superanalytic
$\Lambda_\infty$-superfunction in $(1,0)$-variables}.  If $f(z) \in
\Lambda_\infty^0$ (resp., $\Lambda_\infty^1$) for all $z$ in the
domain of $f$, then $f$ is said to be {\it even} (resp., {\it odd}).
Suppose $f_{(k)} (z_B)$ is convergent in an open neighborhood $N_{(k)}
\subseteq \Bbb{C}$ of $z_B$.  If there exists an open subset $N
\subseteq \bigcap_{(k) \in K_\infty} N_{(k)}$ such that $N \neq
\emptyset$, then $f(z_B)$ is convergent in $N$, and consequently
$f(z)$ is convergent in $\{z = z_B + z_S \in \Lambda_\infty^0 : z_B
\in N \}$.

Let $U \subseteq \Lambda_\infty$ and $H: U \longrightarrow
\Lambda_\infty$, $(z,\theta) \mapsto H(z,\theta)$.  We say that $H$ is
a {\it superanalytic $\Lambda_\infty$- superfunction in
(1,1)-variables} if $H$ is of the form $H(z, \theta) = (f(z) + \theta
\xi(z), \psi(z) + \theta g(z))$ where $f$, $g$, $\xi$, and $\psi$ are
superanalytic $\Lambda_\infty$-superfunctions in $(1,0)$-variables.
If $f$, $g$, $\xi$, and $\psi$ are convergent in the open sets $N_f$,
$N_g$, $N_\xi$, $N_\psi \subseteq \Lambda_\infty^0$, respectively, and
there exists $N_H \subseteq (N_f \cap N_g \cap N_\xi \cap N_\psi)$
such that $N_H \neq \emptyset$, then $H(z,\theta)$ is well-defined
(i.e., convergent) for $\{(z,\theta) \in \Lambda_\infty : z \in N_H
\}$.

Consider the topology on $\Lambda_\infty$ given by the product of the
usual topology on $(\Lambda_\infty)_B = \Bbb{C}$ and the trivial
topology on $(\Lambda_\infty)_S$.  This topology on $\Lambda_\infty$
is called the {\it DeWitt topology}.  The natural domain $U$ of any
superanalytic $\Lambda_\infty$-superfunction is an open set in the
DeWitt topology on $\Lambda_\infty$.

Let $D = \Dz$.  We say that $H(z, \theta) = (\tilde{z},
\tilde{\theta})$ is {\it superconformal} if $D \tilde{z} =
\tilde{\theta} D \tilde{\theta}$.  This is equivalent to 
\[H(z,\theta) = \left(f(z) + \theta \psi(z) \sqrt{f'(z) + \psi(z)
\psi'(z)}, \psi(z) + \theta \sqrt{f'(z) + \psi(z) \psi'(z)} \right)
.\] 
Thus a superconformal $\Lambda_\infty$-superfunction is determined by
an even $\Lambda_\infty$-superfunction in $(1,0)$-variables $f$, an
odd $\Lambda_\infty$-superfunction in $(1,0)$-variables $\psi$, and a
choice of square root for $f_B'(z_B)$.  A square root for $f_B'(z_B)$
is called a {\it spin structure}.  Fix a branch of the complex
logarithm.  For the double valued square root, we will call the
structure given by $\sqrt{1} = 1$ the positive square root structure,
and the structure given by $\sqrt{1} = -1$ the negative square root
structure.  If $H(z,\theta)$ has positive square root structure then
$H(z, - \theta)$ has negative square root structure.  Let
$(\Lambda_\infty^0)^\times = \{a \in \Lambda_\infty^0 : (a_B) \neq 0
\}$.  For $a_0 \in (\Lambda_\infty^0)^\times$, define the linear
operators $a_0^{z \frac{\partial}{\partial z}}$ and $a_0^{ \frac{1}{2}
\theta \frac{\partial}{\partial \theta}}$ on $\Lambda_\infty
[z,z^{-1},\theta]$ by $a_0^{z \frac{\partial}{\partial z}} \cdot c
\theta^m z^n = c \theta^m a_0^n z^n$ and $a_0^{\frac{1}{2} \theta
\frac{\partial}{\partial \theta}} \cdot c \theta^m z^n = c
\sqrt{a_0}^m \theta^m z^n$ for $c \in \Lambda_\infty$, $m \in
\Bbb{Z}_2$, and $n \in \Bbb{Z}$. 

Let $\Lambda_\infty^\infty$ denote the set of infinite series in
$\Lambda_\infty^0 \oplus \Lambda_\infty^1$ indexed by $j \in \Z$ in
the even coordinates and $j - \frac{1}{2}$ for $j \in \Z$ in the odd
coordinates.  We will denote an element of $\Lambda_\infty^\infty$ by
$(A,M) = \{(A_j, M_{j - \frac{1}{2}})\}_{j \in \Z}$ for $A_j \in
\Lambda_\infty^0$ and $M_{j - \frac{1}{2}} \in \Lambda_\infty^1$.  The
following two propositions (proved in \cite{B}) characterize certain
superconformal superfunctions in terms of exponentials of certain
superderivations.  This characterization is analogous to Huang's
\cite{H1} characterization of certain conformal functions in terms of
exponentials of certain derivations.

\begin{prop}\label{at zero} 
Let $(A,M) \in \Lambda_\infty^\infty$, and $a_0 \in
(\Lambda_\infty^0)^\times$.  The formal exponential of formal
differential operators applied to $(z,\theta)$
\begin{equation}\label{superconformal at zero}
H(z, \theta) = \exp\left( \sum_{j \in \Z} \left( A_j \left( \Lz
\right) \right. \right. \hspace{2in}
\end{equation}
\[\hspace{2.5in} \left. \left. + \; M_{j - \frac{1}{2}} \Gz \right)
\right) \cdot a_0^{\left( \Loz \right)} \cdot (z, \theta) \] 
is a formal power series in $z$ and $\theta$ of the form 
\[\left(a_0z + \sum_{j \in \Z} a_j z^{j + 1} + \theta \sum_{j \in \Z}
n_j z^j, \sum_{j \in \Z} m_j z^j + \theta (\sqrt{a_0} + \sum_{j \in
\Z} b_j z^j) \right)\]  
for $a_j, b_j \in \Lambda_\infty^0$ and $m_j, n_j \in
\Lambda_\infty^1$.  If this formal power series converges in a
(DeWitt) open neighborhood in $\Lambda_\infty$, it is superconformal
with positive square root structure, and any superconformal
superfunction vanishing at zero with positive square root structure is
of the form (\ref{superconformal at zero}).  
\end{prop}

\begin{prop}\label{at infinity}
Let $(B,N) = \{(B_j, N_{j - \frac{1}{2}})\}_{j \in \Z}$ be an infinite
series in $\Lambda_\infty^0 \oplus \Lambda_\infty^1$. The formal
exponential of formal differential operators applied to $(\frac{1}{z},
\frac{i \theta}{z})$  
\begin{equation}\label{superconformal at infinity}
H(z, \theta) = \exp\left( - \sum_{j \in \Z} \left( B_j \left(
\Lminusz \right) \right. \right. \hspace{2in}
\end{equation}
\[\hspace{3in} \left. \left. + \; N_{j - \frac{1}{2}} \Gminusz \right)
\right) \cdot \left( \frac{1}{z}, \frac{i \theta}{z} \right)\] 
is a formal power series in $z^{-1}$ and $\theta$ of the form 
\[\left( z^{-1} + \sum_{j \in \Z} a_j z^{- j - 1} + \theta
\sum_{j \in \Z} n_j 
z^{-j - 1}, \sum_{j \in \Z} m_j z^{-j} + \theta (iz^{-1} + \sum_{j \in
\Z} b_j z^{-j - 1}) \right)\]  
for $a_j, b_j \in \Lambda_\infty^0$ and $m_j, n_j \in
\Lambda_\infty^1$.  If this formal power series converges in a
(DeWitt) open neighborhood in $\Lambda_\infty$, it is superconformal
with positive square root structure, and any superconformal
superfunction vanishing at $z_B = \infty$ with positive square root
structure is of the form (\ref{superconformal at infinity}).  
\end{prop}

The two propositions given above can be formulated in terms of the
negative square root structure by the superconformal transformation
$J(z,\theta) = (z,-\theta)$.  

Consider the Neveu-Schwarz Lie superalgebra $\frak{ns}$ (a certain
super-extension of the Virasoro algebra) generated by a central
element $c$, even elements $L(n)$ and odd elements $G(n +
\frac{1}{2})$ for $n \in \Bbb{Z}$ with the following supercommutation 
relations   
\begin{eqnarray*}
[L(m),L(n)] &=& (m - n)L(m + n) + \frac{1}{12} (m^3 - m) \delta_{m + n
, 0} c , \\
\left[ G(m + \frac{1}{2}),L(n) \right] &=& (m - \frac{n - 1}{2} ) G(m
+ n + \frac{1}{2})  \\  
\left[ G (m + \frac{1}{2} ) , G(n - \frac{1}{2} ) \right] &=& 2L(m +
n) + \frac{1}{3} (m^2 + m) \delta_{m + n , 0} c ,
\end{eqnarray*}
for $m, n \in \Bbb{Z}$.

Let $x$ be an {\it even} formal variable (i.e., $x$ commutes with
$\Lambda_\infty$) and $\varphi$ be an {\it odd} formal variable
(i.e., $\varphi$ commutes with $x$ and $\Lambda_\infty^0$,
anti-commutes with $\Lambda_\infty^1$, and $\varphi^2 = 0$).  For any
$s \in \Bbb{C} \diagdown \{0\}$ and $t \in \Bbb{C}$,   
\begin{eqnarray*}
L(n)_t &=& - \left( x^{n + 1} \frac{\partial}{\partial x} + \left(
\frac{n - 1}{2} +  t \right) \varphi x^n \frac{\partial}{\partial
\varphi} \right) \\
G(n + \frac{1}{2})_{t,s} &=& - \left( sx^{n + t}
\frac{\partial}{\partial \varphi} - \frac{1}{s} \varphi x^{n - t +
2}\frac{\partial}{\partial x} \right)  
\end{eqnarray*}
gives a representation of $\frak{ns}$ with $c = 0$.  Propositions
\ref{at zero} and \ref{at infinity} state that formally $L(n)_1$ and
$G(n + \frac{1}{2})_{1,1}$ for $n \in \Bbb{N}$ are the superconformal
infinitesimal transformations at zero with positive square root
structure, and $L(- n)_1$ and $G(- n + \frac{1}{2})_{1,1}$ for $n \in
\Z$ are the superconformal infinitesimal transformations at infinity
with positive square root structure.  The corresponding results for
the negative square root structure state that $L(n)_1$ and $G(n +
\frac{1}{2})_{1,-1} (= - G(n + \frac{1}{2})_{1,1})$ for $n \in
\Bbb{N}$ are the superconformal infinitesimal transformations at zero
with negative square root structure, and $L(- n)_1$ and $G(- n +
\frac{1}{2})_{1,-1} (= - G(- n + \frac{1}{2})_{1,1})$ for $n \in \Z$
are the superconformal infinitesimal transformations at infinity with
negative square root structure.

\section{Superspheres with tubes and the sewing operation}

In this section, we extend Huang's \cite{H1} definition of the moduli
space of spheres with tubes and a sewing operation to a definition of
the moduli space of superspheres with tubes and a sewing operation.
Though it is similar in spirit, we find that in the super case there
is a great deal of non-trivial additional structure involving the soul
coordinates.  

By {\it supersphere} we will mean a supermanifold with DeWitt topology
over $\Lambda_\infty$ (cf. \cite{D}) such that its body is a
genus-zero one-dimensional connected compact complex manifold and its
transition functions are superconformal with a given square root
structure.  A supersphere with $1 + n$ tubes ($n \in \Bbb{N}$) is a
supersphere $S$ with 1 negatively oriented point and $n$ positively
ordered points (called {\it punctures}) and local superconformal
coordinates vanishing at the punctures.  A {\it superconformal
equivalence} $F$ from one supersphere $S_1$ with $1 + n$ tubes to
another supersphere $S_2$ with $1 + n$ tubes which preserves square
root structure is a superconformal isomorphism from the underlying
supersphere of $S_1$ to the underlying supersphere of $S_2$ such that
the $i$-th puncture of $S_1$ is mapped to the $i$-th puncture of
$S_2$, the pull-back of the local coordinate map vanishing at the
$i$-th puncture of $S_2$ is equal to the local coordinate map
vanishing at the $i$-th puncture of $S_1$ in some neighborhood of this
puncture, and the square root structures on $S_1$ and $S_2$ are the
same.

Let $S_1$ be a supersphere with $1 + n$ ($n > 0$) tubes and $S_2$ a
supersphere with $1+ m$ tubes.  Let $p_0,...,p_n$ be the punctures of
$S_1$ with local coordinate charts $(U_i, \Omega_i)$ at $p_i$.  Let
$q_0,...,q_n$ be the punctures of $S_2$ with local coordinate charts
$(V_i, \Xi_i)$ at $q_i$.  We will describe the operation of sewing the
second supersphere at $q_0$ to the first supersphere at $p_i$ for some
fixed $0 \leq i \leq n$, and by the uniformization theorem for
super-Riemann surfaces \cite{CR}, the resulting supermanifold will be
a supersphere.  Assume that there exists a positive number $r$ such
that $\Omega_i(U_i)$ contains the closed set $\bar{\cal B}_0^r =
\bar{B}_0^r \times (\Lambda_\infty)_S$ centered at 0 with radius $r$
in the body, and $\Xi_0 (V_0)$ contains the closed set $\bar{\cal
B}_0^{1/r}$ centered at 0 with radius $1/r$ in the body.  Assume also
that $p_i$ and $q_0$ are the only punctures in $\Omega_i^{-1}
(\bar{\cal B}_0^r)$ and $\Xi_0^{-1} (\bar{\cal B}_0^{1/r})$
respectively.  In this case we say that the $i$-th puncture of the
first supersphere with tubes can be sewn with the $0$-th puncture of
the second supersphere with tubes.  From these two superspheres with
tubes, we can obtain a supersphere with $1 + (n + m - 1)$ tubes by
cutting $\Omega_i^{-1} (\bar{\cal B}_0^r)$ and $\Xi_0^{-1} (\bar{\cal
B}_0^{1/r})$ from $S_1$ and $S_2$ respectively, and then identifying
the boundaries of the resulting surfaces using the map $\Omega_i \circ
I \circ \Xi_0^{-1}$ where $I$ is the map from $\Lambda_\infty^\times$
to itself given by $I(z,\theta) = \left( \frac{1}{z}, \frac{i
\theta}{z} \right)$.  The punctures (with ordering) of this
supersphere with tubes are $p_0, ..., p_{i-1}, q_1, ..., q_m, p_{i +
1}, ..., p_n$.  The local coordinates vanishing at these punctures are
$(U_j, \Omega_j)$ at $p_j$ and $(V_k, \Xi_k)$ at $q_k$.  This
supersphere is denoted $S_1 \; _i\infty_0 \; S_2$ (using Vafa's
\cite{V} notation for the sewing of two spheres).  Note that this
resulting supersphere is independent of the positive number $r$.

The collection of all superconformal equivalence classes of
superspheres with tubes and a given square root structure is called
{\it the moduli space of superspheres with tubes}.  The global
superconformal transformations with positive (resp., negative) square
root structure are generated by $L(\pm 1)_1$, $L(0)_1$, and
$G(\pm\frac{1}{2})_{1,1}$ (resp., $G(\pm\frac{1}{2})_{1,-1} = -
G(\pm\frac{1}{2})_{1,1}$).  Note that since in defining the sewing
operation of two superspheres, we chose the identification of
boundaries to be $I(z,\theta) = \left( \frac{1}{z}, \frac{i \theta}{z}
\right)$, the sewing of two superspheres with the same square root
structure results in a supersphere with the same square root structure
as the original two superspheres.  (If we had chosen $I^{-1}(z,\theta)
= \left( \frac{1}{z}, -\frac{i \theta}{z} \right)$, the sewing of two
superspheres with the same square root structure would result in a
supersphere with opposite square root structure.)  Thus the sewing
operation described above is a well-defined (partial) operation on the
moduli space of superspheres with tubes and a given square root
structure.  The operation is partial since not all superspheres with
tubes can be sewn together.  As in \cite{H1} and \cite{H2}, this is an
important feature in that the non-sewability of certain spheres
reflects information related to the analytic structure of the the
moduli space.  For a given branch cut in the complex plane, there are
two square root structures, and thus we obtain two moduli spaces --
one with a positive square root structure and one with a negative
square root structure with respect to the branch cut.  We can
explicitly describe these moduli spaces using Propositions \ref{at
zero} and \ref{at infinity} to describe the local coordinates.  Let
$(\Lambda_\infty^0)^\times \times {\cal H}$ be the subset of
$(\Lambda_\infty^0)^\times \times \Lambda_\infty^\infty$ consisting of
all elements $(a_0,A,M)$ such that (\ref{superconformal at zero}) is a
convergent power series in some neighborhood of zero, ${\cal H}^{(0)}$
the subset of $\Lambda_\infty^\infty$ containing all elements $(B,N)$
such that (\ref{superconformal at infinity}) is a convergent power
series in some neighborhood of infinity, and $SM^{n - 1}$ the subset
of elements in $\Lambda_\infty^{n-1}$ with distinct non-zero bodies.
In \cite{B}, it is shown that the moduli space of superspheres with $1
+ n$ tubes ($n > 0$) and a given square root structure can be
identified with the set $SK(n) = SM^{n-1} \times {\cal H}^{(0)} \times
((\Lambda_\infty^0)^\times \times {\cal H})^n$, and the moduli space
of superspheres with one tube and a given square root structure can be
identified with the set $SK(0) = \{(B,N) \in {\cal H}^{(0)} : (A_1,
M_{\frac{1}{2}}) = (0,0) \}$.  Thus we will refer to $SK^+$
(resp., $SK^-$) as the moduli space of superspheres with tubes and
positive (resp., negative) square root structure.  Although as sets
$SK^+ = SK^-$, the underlying geometric structures they represent are
distinct, and we can map one to the other via the transformation
$J(z,\theta) = (z,-\theta)$ on the underlying superspheres with tubes.
The definition of sewing of superspheres gives a (partial) operation
on $SK^+$ (resp., $SK^-$) which we again denote by $_i\infty_0$, and as
sets, the sewing operation on $SK^+$ is the same as that on $SK^-$.
We write an element of $SK^+(n)$ or $SK^-(n)$ as $((z_1,
\theta_1),...,(z_{n-1}, \theta_{n-1}); (A^{(0)},M^{(0)}), (a_0^{(1)},
A^{(1)}, M^{(1)}),...,(a_0^{(n)}, A^{(n)}, M^{(n)}))$.  (When we want
to refer to $SK^+$ or $SK^-$, we will write $SK^\pm(n)$.)  The
symmetric group on $n - 1$ letters $S_{n-1}$ acts on $SK^\pm(n)$ by
permuting the $(z_i, \theta_i)$ and $(a_0^{(i)}, A^{(i)}, M^{(i)})$
for $i = 1,...,n-1$.  In \cite{B}, we extend this action to an action
of $S_n$ on $SK^\pm(n)$.  We will denote by {\bf 0} the element of
$\Lambda_\infty^\infty$ with all components equal to 0.  The
superconformal transformation $J(z,\theta) = (z, -\theta)$ induces a
map (which we also denote by $J$) from $SK^\pm$ to $SK^\mp$ given by
\[J((z_1,\theta_1),...,(z_{n-1}, \theta_{n-1}); (A^{(0)},M^{(0)}),
(a_0^{(1)}, A^{(1)}, M^{(1)}),...,(a_0^{(n)}, A^{(n)}, M^{(n)})) =
\hspace{.7in} \]
\[\hspace{.9in}
((z_1,-\theta_1),...,(z_{n-1},-\theta_{n-1});(A^{(0)},M^{(0)}),
(a_0^{(1)}, A^{(1)}, M^{(1)}),...,(a_0^{(n)}, A^{(n)}, M^{(n)})) . \] 

\begin{rema} The moduli space $SK^\pm$ of superspheres with tubes
and a given positive or negative square root structure and the sewing
operation on $SK^\pm$  is a partial operad (cf. \cite{HL}).
\end{rema}

Let $Q_1 \in SK^\pm (n)$ for $n \in \Z$, and $Q_2 \in SK^\pm (m)$ for
$m \in \Bbb{N}$.  Let the coordinates at the puncture $(z_i,
\theta_i)$ of $Q_1$ be $H(z - z_i - \theta \theta_i, \theta -
\theta_i)$ where $H(z,\theta)$ is given by (\ref{superconformal at
zero}), and let the local coordinates at infinity of $Q_2$ be given by
(\ref{superconformal at infinity}).  The following proposition
describes the change of local coordinates of the resulting supersphere
$Q_1 \; _i\infty_0 \; Q_2$, and this description is given in terms of
elements of the Neveu-Schwarz algebra with central charge $c \in
\Bbb{C}$. 

\begin{prop}\label{Gamma and Psi}  
Let $({\cal A}, {\cal M}) = \{ ({\cal A}_j, {\cal M}_{j -
\frac{1}{2}}) \}_{j \in \Z}$ and $({\cal B}, {\cal N}) = \{ ({\cal
B}_j, {\cal N}_{j - \frac{1}{2}}) \}_{j \in \Z}$ be two sequences of
formal variables, ${\cal A}_j$ and ${\cal B}_j$ even and ${\cal M}_{j
- \frac{1}{2}}$ and ${\cal N}_{j - \frac{1}{2}}$ odd, let $\alpha_0$
be another even formal variable, and let $V$ be a positive energy
module for the Neveu-Schwarz algebra.  There exist unique 
canonical series $(\Psi_j,\Psi_{j - \frac{1}{2}}) = (\Psi_j,\Psi_{j -
\frac{1}{2}})(\alpha_0, {\cal A}, {\cal M}, {\cal B}, {\cal N})$ for
$j \in \Bbb{Z}$, and $\Gamma = \Gamma (\alpha_0, {\cal A}, {\cal M},
{\cal B}, {\cal N})$ in $\Bbb{C} [\alpha_0, \alpha_0^{-1},
\sqrt{\alpha_0}, \sqrt{\alpha_0}^{-1}] [[{\cal A}, {\cal M}, {\cal B},
{\cal N}]]$ such that

\[e^{- \sum_{j \in \Z} ({\cal A}_j L(j) \; + \; {\cal M}_{j -
\frac{1}{2}} G(j - \frac{1}{2}))} \; \alpha_0^{-L(0)} \; e^{- \sum_{j
\in \Z} ({\cal B}_j L(-j) \; + \; {\cal N}_{j - \frac{1}{2}} G(- j +
\frac{1}{2}))} = \hspace{1.4in} \] 
\[\hspace{.8in} e^{\sum_{j \in \Z} (\Psi_{-j} L(-j) \; + \; \Psi_{- j
+ \frac{1}{2}} G(- j + \frac{1}{2}))} \; e^{\sum_{j \in \Z} (\Psi_j
L(j) \; + \; \Psi_{j - \frac{1}{2}} G(j - \frac{1}{2}))} \; e^{\Psi_0
L(0)} \; \alpha_0^{-L(0)} \; e^{\Gamma c} \]  
as operators in $(\mbox{\em End} \; V) [\alpha_0, \alpha_0^{-1},
\sqrt{\alpha_0}, \sqrt{\alpha_0}^{-1}] [[{\cal A}, {\cal M}, {\cal B},
{\cal N}]]$.
\end{prop}

The series $\Gamma$ can easily be calculated up to second order terms
in the ${\cal A}_j$'s, ${\cal M}_{j - \frac{1}{2}}$'s, ${\cal B}_j$'s,
and ${\cal N}_{j - \frac{1}{2}}$'s for $j \in \Z$.  In fact,
\begin{eqnarray*}
\Gamma &=& \Gamma(\alpha_0,{\cal A},{\cal M},{\cal B},{\cal N})  \\  
&=& \sum_{j \in \Z} \left( \left( \frac{j^3 - j}{12} \right)
\alpha_0^{-j} {\cal A}_j {\cal B}_j + \left( \frac{j^2 - j}{3} \right)
\sqrt{\alpha_0} \alpha_0^{-j} {\cal N}_{j - \frac{1}{2}} {\cal M}_{j -
\frac{1}{2}} \right) + \Gamma_0
\end{eqnarray*}
where $\Gamma_0$ contains only terms with products of at least three
of the ${\cal A}_j$'s, ${\cal M}_{j - \frac{1}{2}}$'s, ${\cal B}_j$'s,
and ${\cal N}_{j - \frac{1}{2}}$'s for $j \in \Z$ (but not all of the
three ${\cal A}_j$'s, ${\cal M}_{j - \frac{1}{2}}$'s, ${\cal B}_j$'s,
or ${\cal N}_{j - \frac{1}{2}}$'s).  The series $\Gamma$ has the
following convergence property:

\begin{prop}\label{convergence of t-series}
Let
\[Q_1^\pm = ((z_1, \theta_1),...,(z_{n-1}, \theta_{n-1});
(A^{(0)},M^{(0)}), (a_0^{(1)}, A^{(1)}, M^{(1)}),...,(a_0^{(n)},
A^{(n)}, M^{(n)})) \in SK^\pm(n), \]
\[Q_2^\pm = ((z_1, \theta_1),...,(z_{m-1}, \theta_{m-1});
(B^{(0)},N^{(0)}), (b_0^{(1)}, B^{(1)}, N^{(1)}),...,(b_0^{(m)},
B^{(m)}, N^{(m)})) \! \in SK^\pm(m) .\]
If the i-th tube of $Q_1^\pm$ can be sewn with the 0-th tube of
$Q_2^\pm$, the t-series
\[e^{\Gamma (t^{-1}a_0^{(i)}, A^{(i)}, M^{(i)}, B^{(0)}, N^{(0)})c}\]
is absolutely convergent at $t = 1$.
\end{prop}

\section{The linear algebra of $\mathbf{Z}_2$-graded
$\Lambda_\infty$-modules with $\frac{1}{2} \mathbf{Z}$-graded
finite-dimensional weight spaces}

For any $\Bbb{Z}_2$-graded vector space $V = V^0 \oplus V^1$, define
the {\it sign} of $v$ homogeneous in $V$ to be $\eta(v) = i$ for $v
\in V^i$, $i \in \Bbb{Z}_2$.  Let
\[V = \coprod_{n \in \frac{1}{2} \Bbb{Z}} V_{(n)} = \coprod_{n \in
\frac{1}{2} \Bbb{Z}} V_{(n)}^0 \oplus \coprod_{n \in \frac{1}{2}
\Bbb{Z}} V_{(n)}^1 = V^0 \oplus V^1 \]  
with
\[\dim V_{(n)} < \infty \quad \mbox{for} \quad n \in \frac{1}{2}
\Bbb{Z} ,\] 
be a $\frac{1}{2} \Bbb{Z}$-graded (by weight) $\Lambda_\infty$-module
with finite-dimensional homogeneous weight spaces $V_{(n)}$ which is
also $\Bbb{Z}_2$-graded (by sign).  Let 
\[V' = \coprod_{n \in \frac{1}{2} \Bbb{Z}} V_{(n)}^* \]
be the graded dual space of $V$, 
\[\bar{V} = \prod _{n \in \frac{1}{2} \Bbb{Z}} V_{(n)} = V'^* \] 
the algebraic completion of $V$, and $\langle \cdot , \cdot \rangle$
the natural pairing between $V'$ and $\bar{V}$.  For any $n \in
\Bbb{N}$, let 
\[{\cal SF}_V(n) = \mbox{Hom}_{\Lambda_\infty}(V^{\otimes n}, \bar{V})
.\] 
For any $m \in \Z$, $n \in \Bbb{N}$, and any positive integer $i \leq
m$, we define the {\it $t$-contraction} 
\begin{eqnarray*} 
 _i*_0 : {\cal SF}_V(m) \times {\cal SF}_V(n) & \rightarrow &
\mbox{Hom} (V^{\otimes (m + n - 1)}, V[[t^{\frac{1}{2}}, t^{-
\frac{1}{2}}]])  \\ 
(f,g) &\mapsto& (f \; _i*_0 \; g)_t , 
\end{eqnarray*}
by
\[(f \; _i*_0 \; g)_t(v_1 \otimes \cdots \otimes v_{n + m - 1}) =
\hspace{4in} \]
\begin{equation}\label{t-contraction}
\sum_{k \in \frac{1}{2} \Bbb{Z}}f(v_1 \otimes \cdots \otimes
v_{i - 1} \otimes P_k(g(v_i \otimes \cdots \otimes v_{i + n - 1}))
\otimes v_{i + n} \otimes \cdots \otimes v_{m + n - 1}) t^k 
\end{equation}
for all $v_i,...,v_{m + n - 1} \in V$, where for any $k \in
\frac{1}{2} \Bbb{Z}$, $P_k : V \rightarrow V_{(k)}$ is the projection
map.  If we want to substitute complex values for $t$ into equation
(\ref{t-contraction}), we must choose a square root.  However, we have
already fixed a branch cut and defined the positive and negative
single-valued square roots for this branch cut.  We denote the two
corresponding positive and negative $t$-contractions by $(f \; _i*_0
\; g)^+_t$ and $(f \; _i*_0 \; g)^-_t$, respectively, for $t \in
\Bbb{C}$.    

If for arbitrary $v' \in V'$, $v_1,...,v_{m + n -1} \in V$, the formal
Laurent series in $t^{\frac{1}{2}}$
\[ \langle v', (f \; _i*_0 \; g)^\pm_t(v_1 \otimes \cdots \otimes v_{n
+ m - 1}) \rangle \]
is absolutely convergent when $t = 1$, then $(f \; _i*_0 \; g)^\pm_1$
is well-defined as an element of ${\cal SF}_V(m + n -1)$, and we
define the {\it positive (resp., negative) contraction $(f \; _i*_0 \;
g)^+$, (resp. $(f \; _i*_0 \; g)^-$) in ${\cal SF}_V(m + n -1)$ of $f$
and $g$} by  
\[(f \; _i*_0 \; g)^\pm = (f \;  _i*_0 \; g)^\pm_1 .\]

Let $(l \; k) \in S_n$ be the permutation on $n$ letters which
switches the $l$-th and $k$-th letters, for $l, k = 1, ..., n$, $l <
k$.   We define an action of the transposition $(l \; k)$ on
$V^{\otimes n}$ by    
\[ (l \; k)(v_1 \otimes \cdots \otimes v_l \otimes \cdots \otimes v_k
\otimes \cdots  \otimes v_n) = (-1)^{\eta{(l \; k)}} (v_1 \otimes
\cdots \otimes v_k \otimes \cdots \otimes v_l \otimes \cdots \otimes
v_n) \]
for $v_j$ of homogeneous sign in $V$, where 
\[\eta{(l \; k)} = \sum_{j = l + 1}^{k-1} \eta(v_j)(\eta(v_l) +
\eta(v_k)) + \eta(v_k)\eta(v_l) .\]
Let $\sigma \in S_n$ be a permutation on $n$ letters. Then $\sigma$ is
the product of transpositions $\sigma = \sigma_1 \cdots 
\sigma_m$, $\sigma_i = (l_i \; k_i)$, $l_i, k_i \in \{1,...,n\}$,
$l_i < k_i$, $i = 1,...,m$.   Thus we have an action of $S_n$ on
$V^{\otimes n}$ given by  
\[\sigma(v_1 \otimes \cdots \otimes v_n) = \sigma_1 \cdots
\sigma_m (v_1 \otimes \cdots \otimes v_n) = (-1)^{\eta(\sigma_1) +
\cdots \eta(\sigma_m)} v_{\sigma(1)} \otimes \cdots \otimes
v_{\sigma(n)} .\] 
This action of $S_n$ induces a left action of $S_n$ on ${\cal
SF}_V(n)$ given by 
\[\sigma (f) (v_1 \otimes \cdots \otimes v_n) = f(\sigma^{-1}(v_1
\otimes \cdots \otimes v_n)) , \]
for $f \in {\cal SF}_V(n)$.

Since $V$ is a $\Bbb{Z}_2$-graded $\Lambda_\infty$-module, End $V$ has
a natural $\Bbb{Z}_2$-grading given by even operators $(\mbox{End} \;
V)^0 = \{P \in \mbox{End} \; V : PV^m \subset V^m \; \mbox{for} \; m
\in \Bbb{Z}_2 \}$ and odd operators $(\mbox{End} \; V)^1 = \{P \in
\mbox{End} \; V : PV^m \subset V^{(m + 1) \mbox{\begin{footnotesize}
mod \end{footnotesize}} 2} \; \mbox{for}
\; m \in \Bbb{Z}_2 \}$.  Also, End $V$ has a natural Lie superalgebra
structure with supercommutator given by $[P_1, P_2] = P_1 P_2 -
(-1)^{\eta(P_1) \eta(P_2)} P_2 P_1$, for $P_1$ and $P_2$ of homogeneous
sign in End $V$. 

If $P \in \mbox{End} \; V$, the corresponding adjoint operator on
$V'$, if it exists, is denoted by $P'$.  The condition for the
existence of $P'$ is that the linear functional on $V$ defined by the 
right-hand side of  
\[ \langle P'v',v\rangle = \langle v',Pv \rangle, \; \mbox{for} \; v
\in V, \; v' \in V' \]
should lie in $V'$.  If there exists $n \in \frac{1}{2} \Bbb{Z}$ such
that $P$ maps $V_{(k)}$ to $V_{(n + k)}$ for any $k \in \frac{1}{2}
\Bbb{Z}$, we say that $P$ has {\it weight $n$}.  It is easy to see
that $P$ has weight $n$ if and only if its adjoint $P'$ exists and has
weight $-n$ as an operator on $V'$, and that $P$ is even (resp., odd)
if and only if its adjoint exists and is even (resp., odd).  In the
case that $V$ is a module for the Neveu-Schwarz algebra graded by the 
eigenvalues of $L(0)$, the adjoint operator $L'(-n)$ for $n \in
\Bbb{Z}$ corresponding to $L(-n)$ exists and is even with weight $n$,
and the adjoint operator $G'(-n - \frac{1}{2})$ for $n \in \Bbb{Z}$
corresponding to $G(-n - \frac{1}{2})$ exists and is odd with weight
$n + \frac{1}{2}$.

\section{Supergeometric vertex operator superalgebras}

In the definition of supergeometric vertex operator superalgebras, we
need the following notion of supermeromorphic superfunction on
$SK^\pm(n)$.  A {\it supermeromorphic superfunction on $SK^\pm(n)$}
($n \in \Z$) is a superfunction $F : SK^\pm(n) \rightarrow
\Lambda_\infty$ of the form 
\[F((z_1, \theta_1),...,(z_{n-1}, \theta_{n-1}); (A^{(0)},M^{(0)}),
(a_0^{(1)}, A^{(1)}, M^{(1)}),...,(a_0^{(n)}, A^{(n)}, M^{(n)}))
\hspace{1in} \]
\begin{equation}\label{supermeromorphic}
= \frac{1}{\prod^{n-1}_{i = 1} z_i^{s_i} \prod_{1 \leq i<j \leq n-1}
(z_i - z_j - \theta_i \theta_j)^{s_{ij}}} F_0((z_1,
\theta_1),...,(z_{n-1}, \theta_{n-1}); (A^{(0)},M^{(0)}),
\end{equation}
\[\hspace{3.5in} (a_0^{(1)}, A^{(1)}, M^{(1)}),...,(a_0^{(n)},A^{(n)},
M^{(n)}))\] 
where $s_i$ and $s_{ij}$ are nonnegative integers and 
\[F_0((z_1,\theta_1),...,(z_{n-1}, \theta_{n-1}); (A^{(0)},M^{(0)}),
(a_0^{(1)}, A^{(1)}, M^{(1)}),...,(a_0^{(n)},A^{(n)}, M^{(n)}))\]  
is a polynomial in the $z_i$'s, $\theta_i$'s, $a_0^{(i)}$'s,
$(a_0^{(i)})^{-1}$'s, $A^{(i)}_j$'s, and $M^{(i)}_{j -
\frac{1}{2}}$'s.  For $n=0$ a {\it supermeromorphic superfunction on
$SK^\pm(0)$} is a polynomial in the components of elements of
$SK^\pm(0)$.   

For $L \in \Bbb{N}$, let $\Lambda_L$, be the Grassmann subalgebra over
{\bbb C} of $\Lambda_\infty$ on generators $\zeta_1, \zeta_2,
...,\zeta_L$.  We use the notation $\Lambda_*$ to denote
$\Lambda_L$ for some $L \in \Bbb{N}$ or $\Lambda_\infty$.  Recall that
a formal variable is even if it commutes with $\Lambda_\infty$ and all
other formal variables and is odd if it commutes with
$\Lambda_\infty^0$, and anti-commutes with $\Lambda_\infty^1$ and all
odd formal variables including itself, i.e., its square is zero.

\begin{defn} {\em A} ($N = 1$ Neveu-Schwarz) supergeometric vertex
operator superalgebra over $\Lambda_*$ with positive square root
structure {\em is a $\frac{1}{2} \Bbb{Z}$-graded (by weight)
$\Lambda_\infty$-module which is also $\Bbb{Z}_2$-graded (by sign)  
 \[V = \coprod_{n \in \frac{1}{2} \Bbb{Z}} V_{(n)} = \coprod_{n \in
 \frac{1}{2} \Bbb{Z}} V_{(n)}^0 \oplus \coprod_{n \in  \frac{1}{2}
\Bbb{Z}} V_{(n)}^1 = V^0 \oplus V^1 \]  
such that only the subspace $\Lambda_*$ of $\Lambda_\infty$ acts
non-trivially on $V$, 
\[\dim V_{(n)} < \infty \quad \mbox{for} \quad n \in \frac{1}{2}
\Bbb{Z} ,\] 
and for any $n \in \Bbb{N}$, a map
\[ \nu_n^+ : SK^+(n) \rightarrow S{\cal F}_V (n) \]
satisfying the following axioms: 

(1)} Positive energy axiom: 
{\em \[V_{(n)} = 0 \quad \mbox{for $n$ sufficiently small.} \]

(2)} Grading axiom:  {\em Let $v' \in V'$, $v \in V_{(n)}$, and $a_0
\in (\Lambda_\infty^0)^\times$.  Then
\[ \langle v', \nu_1^+(\mbox{\bf 0},(a_0,\mbox{\bf 0}))(v)\rangle =
a_0^{-n} \langle v',v \rangle .\]

(3)} Supermeromorphicity axiom: {\em For any $n \in \Z$, $v' \in V'$,
and $v_1,...,v_n \in V$, the function
\[Q \mapsto \langle v', \nu_n^+(Q) (v_1 \otimes \cdots \otimes v_n)
\rangle \]
on $SK^+(n)$ is a canonical supermeromorphic superfunction (in the
sense of (\ref{supermeromorphic})), and if $(z_i, \theta_i)$ and
$(z_j, \theta_j)$ are the i-th and j-th punctures of $Q \in SK^+(n)$
respectively ($i,j \in \{1,...,n\}$, $i \neq j$), then for any $v_i$
and $v_j$ in $V$ there exists $N(v_i,v_j) \in \Z$ such that for any
$v' \in V'$ and $v_k \in V$, $k \neq i,j$, the order of the pole
$(z_i, \theta_i) = (z_j, \theta_j)$ of $\langle v', \nu_n(Q) (v_1
\otimes \cdots \otimes v_n) \rangle$ is less then $N(v_i,v_j)$.  

(4)} Permutation axiom:  {\em Let $\sigma \in S_n$.  Then for any $Q
\in SK^+(n)$
\[\sigma(\nu_n^+ (Q)) = \nu_n^+(\sigma(Q)) . \]

(5)} Sewing axiom:  {\em There exists a unique complex number $c$ (the}
central charge {\em or} rank{\em ) such that if 
\[Q_1 = ((z_1, \theta_1),...,(z_{n-1}, \theta_{n-1});
(A^{(0)},M^{(0)}), (a_0^{(1)}, A^{(1)}, M^{(1)}),...,(a_0^{(n)},
A^{(n)}, M^{(n)})) \in SK^+(n), \]
\[Q_2 = ((z_1', \theta_1'),...,(z_{m-1}', \theta_{m-1}');
(B^{(0)},N^{(0)}), (b_0^{(1)}, B^{(1)}, N^{(1)}),...,(b_0^{(m)},
B^{(m)}, N^{(m)})) \in SK^+(m) ,\] and if the i-th tube of $Q_1$ ($1
\leq i \leq n$) can be sewn with the 0-th tube of $Q_2$, then for any
$v' \in V'$, $v_1,...,v_{n + m - 1} \in V$,
\[\langle v', (\nu_n^+ (Q_1) \; _i*_0 \; \nu_m^+ (Q_2))^+_t (v_1 \otimes 
\cdots \otimes v_{n + m - 1}) \rangle \]
is absolutely convergent when $t = 1$, and
\[\nu_{n + m - 1}^+ (Q_1 \; _i\infty_0 \; Q_2) = (\nu_n^+ (Q_1) \;
_i*_0  \; \nu_m^+ (Q_m))^+ \; e^{-\Gamma (a_0^{(i)}, A^{(i)}, M^{(i)},
B^{(0)}, N^{(0)})c} . \]}
\end{defn}

We denote the supergeometric vertex operator superalgebra defined
above by $(V, \nu^+ = \{\nu_n^+ \}_{n \in \Bbb{N}})$. Replacing $SK^+$
by $SK^-$ in the above definition, we have the corresponding notion of
a supergeometric vertex operator superalgebra with negative square
root structure $(V,\nu^-)$ where $\nu^- = \nu^+ \circ J$.  
 
Let $(V_1, \nu^\pm)$ and $(V_2, \mu^\pm)$ be two supergeometric vertex
operator superalgebras over $\Lambda_*$, $\gamma : V_1 \rightarrow V_2$
a doubly graded linear homomorphism, i.e., $\gamma : (V_1)_{(n)}^i
\rightarrow  (V_2)_{(n)}^i$ for $n \in \frac{1}{2} \Bbb{Z}$, and $i
\in \Bbb{Z}_2$, and let $\bar{\gamma} : \bar{V_1} \rightarrow
\bar{V_2}$ be the unique extension of $\gamma$.  If for any $n \in
\Bbb{N}$ and any $Q \in SK^\pm(n)$
\[\bar{\gamma} \circ \nu^\pm (Q) = \mu^\pm (Q) \circ \gamma^{\otimes
n} , \] 
we say that $\gamma$ is a {\it homomorphism} from  $(V_1, \nu^\pm)$ to
$(V_2, \mu^\pm)$.

Let $c$ be a complex number, and let $\bf{SG}^+(c,*)$ be the category
of supergeometric vertex operator superalgebras over $\Lambda_*$ with
positive square root structure, and $\bf{SG}^-(c,*)$ be the category of
supergeometric vertex operator superalgebras over $\Lambda_*$ with
negative square root structure.  
  
\begin{prop}\label{supergeometric isomorphism}
The two categories $\bf{SG}^+(c,*)$ and $\bf{SG}^-(c,*)$ are
isomorphic. 
\end{prop}

{\it Proof}:  We define $J^+ : \bf{SG}^+(c,*) \rightarrow
\bf{SG}^-(c,*)$ and $J^- : \bf{SG}^-(c,*) \rightarrow  \bf{SG}^+(c,*)$
by  
\begin{eqnarray*}
J^+(V,\nu^+) = (V,\nu^- = \nu^+ \circ J) \quad \mbox{and} \quad
J^+(\gamma) = \gamma \\ 
J^-(V,\nu^-) = (V,\nu^+ = \nu^- \circ J) \quad \mbox{and} \quad
J^-(\gamma) = \gamma . 
\end{eqnarray*}
It is easy to see that $J^+$ and $J^-$ are functors and that $J^+
\circ J^- = 1_{\bf{SG}^-(c,*)}$ and $J^- \circ J^+ =
1_{\bf{SG}^+(c,*)}$.  $\pfbox$

\section{(Superalgebraic) vertex operator superalgebras}

In this section, we extend the notion of vertex operator superalgebra
to include odd formal variables.  In \cite{FLM}, the formal
$\delta$-function $\delta(x) = \sum_{n \in \Bbb{Z}} x^n$ is 
a fundamental ingredient in the formal calculus underlying the theory
of vertex operator algebras.  Extending the formal calculus to include
odd formal variables, we note that for any formal Laurent series $f(x)
\in \Lambda_\infty [[x,x^{-1}]]$ in the even formal variable $x$, and
for any odd formal variables $\varphi_1$ and $\varphi_2$, we have $f(x
+ \varphi_1 \varphi_2) = f(x) + \varphi_1 \varphi_2 f'(x)$.  Thus we
have the following $\delta$-function involving three even variables
and two odd variables: 
\begin{eqnarray*}
\delta \left( \frac{x_1 - x_2 - \varphi_1 \varphi_2}{x_0} \right) &=&
\sum_{n \in \Bbb{Z}} \frac{(x_1 - x_2 - \varphi_1 \varphi_2)^n}{x_0^n}
= \sum_{n \in \Bbb{Z}} \left( \frac{(x_1 - x_2)^n}{x_0^n} - n
\varphi_1 \varphi_2 \frac{(x_1 - x_2)^{n - 1}}{x_0^n} \right) \\
&=& \delta \left( \frac{x_1 - x_2}{x_0} \right)  - \varphi_1 \varphi_2
x_0^{-1} \delta' \left( \frac{x_1 - x_2}{x_0} \right).
\end{eqnarray*}

\begin{defn} {\em A} ($N = 1$ Neveu-Schwarz) vertex operator
superalgebra over $\Lambda_*$ and with odd variables {\em is a
$\frac{1}{2}${\bbb Z}-graded (by weight) $\Lambda_\infty$-module which
is also $\Bbb{Z}_2$-graded (by sign)  
\[V = \coprod_{n \in \frac{1}{2} \Bbb{Z}} V_{(n)} = \coprod_{n
\in \frac{1}{2}\Bbb{Z}} V_{(n)}^0 \oplus \coprod_{n \in
\frac{1}{2}\Bbb{Z}} V_{(n)}^1 = V^0 \oplus V^1 \]  
such that only the subspace $\Lambda_*$ of $\Lambda_\infty$ acts
non-trivially on $V$, and
\[\dim V_{(n)} < \infty \quad \mbox{for} \quad n \in \frac{1}{2}
\Bbb{Z} , \]
\[V_{(n)} = 0 \quad \mbox{for $n$ sufficiently small} , \]
equipped with a linear map $V \otimes V \rightarrow V[[x,x^{-1}]]
\oplus \varphi V[[x,x^{-1}]]$, or equivalently,
\begin{eqnarray*} 
V &\rightarrow&  (\mbox{End} \; V)[[x,x^{-1}]] \oplus \varphi
(\mbox{End} \; V)[[x,x^{-1}]] \\
v  &\mapsto&  Y(v,(x,\varphi)) = \sum_{n \in \Bbb{Z}} v_n x^{-n-1} +
\varphi \sum_{n \in \Bbb{Z}} v_{n - \frac{1}{2}} x^{-n-1}
\end{eqnarray*}
where $v_n \in (\mbox{End} \; V)^{\eta(v)}$ and $v_{n -
\frac{1}{2}} \in (\mbox{End} \; V)^{(\eta(v) + 1)
\mbox{\begin{footnotesize} mod \end{footnotesize}} 2}$ for $v$ of 
homogeneous sign in $V$, $x$ is an even formal variable, and $\varphi$
is an odd formal variable, and where $Y(v,(x,\varphi))$ denotes the}
vertex operator associated with $v$, {\em and equipped also with two
distinguished homogeneous vectors $\mbox{\bf 1} \in V_{(0)}^0$ (the
{\em vacuum}) and $\tau \in V_{(\frac{3}{2})}^1$.  The following
conditions are assumed for $u,v \in V$:   
\[u_n v = 0 \quad \mbox{for $n \in \frac{1}{2} \Bbb{Z}$ sufficiently
large;} \]
\[Y(\mbox{\bf 1}, (x, \varphi)) = 1 \quad \mbox{(1 on the right being
the identity operator);} \]
the} creation property {\em holds:
\[Y(v,(x,\varphi)) \mbox{\bf 1} \in V[[x]] \oplus \varphi V[[x]]
\qquad \mbox{and} \qquad \lim_{(x,\varphi) \rightarrow 0}
Y(v,(x,\varphi)) \mbox{\bf 1} = v ; \]
the} Jacobi identity {\em holds: 
\[x_0^{-1} \delta \left( \frac{x_1 - x_2 - \varphi_1 \varphi_2}{x_0}
\right) Y(u,(x_1, \varphi_1))Y(v,(x_2, \varphi_2)) \hspace{2.8in} \]
\[- (-1)^{\eta(u)\eta(v)} x_0^{-1} \delta \left( \frac{x_2 - x_1 + 
\varphi_1 \varphi_2}{-x_0} \right)Y(v,(x_2, \varphi_2))Y(u,(x_1,
\varphi_1)) \]
\[\hspace{2in} = x_2^{-1} \delta \left( \frac{x_1 - x_0 - \varphi_1
\varphi_2}{x_2} \right) Y(Y(u,(x_0, \varphi_1 - \varphi_2))v,(x_2,
\varphi_2)) ; \]
the Neveu-Schwarz algebra relations hold:
\begin{eqnarray*}
\left[L(m),L(n) \right] &=& (m - n)L(m + n) + \frac{1}{12} (m^3 - m)
\delta_{m + n , 0} c , \\ \label{V9}
\left[ G(m + \frac{1}{2}),L(n) \right] &=& (m - \frac{n - 1}{2} ) G(m
+ n + \frac{1}{2}) ,\\ \label{V10}
\left[ G (m + \frac{1}{2} ) , G(n - \frac{1}{2} ) \right] &=& 2L(m +
n) + \frac{1}{3} (m^2 + m) \delta_{m + n , 0} c , \label{V11}
\end{eqnarray*}
for $m,n \in \Bbb{Z}$, where 
\[G(n + \frac{1}{2}) = \tau_{n + 1}, \qquad \mbox{and} \qquad 2L(n) =
\tau_{n + \frac{1}{2}} \qquad \mbox{for} \; n \in \Bbb{Z} , \]
i.e. 
\[Y(\tau,(x,\varphi)) = \sum_{n \in \Bbb{Z}} G (n + \frac{1}{2}) x^{-
n - \frac{1}{2} - \frac{3}{2}} \; + \; 2 \varphi \sum_{n \in \Bbb{Z}}
L(n) x^{- n - 2} ,\]
and $c \in \Bbb{C}$; 
\[L(0)v = nv \quad \mbox{for} \quad n \in \frac{1}{2} \Bbb{Z} \quad
\mbox{and} \quad v \in V_{(n)}; \]
\[\left( \frac{\partial}{\partial \varphi} + \varphi
\frac{\partial}{\partial x} \right) Y(v,(x,\varphi)) =  Y(G(- 
\frac{1}{2})v,(x,\varphi)) . \] }
\end{defn}
The superalgebraic vertex operator superalgebra just defined is
denoted by $(V,Y(\cdot,(x,\varphi)),\mbox{\bf 1},\tau)$.  

Two consequences of the definition are that 
\begin{equation}\label{Y with x and phi}
Y(v,(x,\varphi)) = \sum_{n \in \Bbb{Z}} v_n x^{-n-1} + \varphi
\sum_{n \in \Bbb{Z}} [G(- \frac{1}{2}), v_n] x^{-n-1} ,
\end{equation}
i.e. $v_{n - \frac{1}{2}} =  [G(- \frac{1}{2}), v_n]$, and
\begin{equation}\label{L derivative}
\frac{\partial}{\partial x} Y(v, (x,\varphi)) =
Y(L(-1)v,(x,\varphi)) .
\end{equation}

Let $(V_1, Y_1(\cdot,(x,\varphi)),\mbox{\bf 1}_1,\tau_1)$ and $(V_2,
Y_2(\cdot,(x,\varphi)),\mbox{\bf 1}_2,\tau_2)$ be two vertex operator
superalgebras over $\Lambda_*$.  A {\it homomorphism} of vertex
operator superalgebras with odd formal variables is a doubly graded
$\Lambda_*$-module homomorphism $\gamma : V_1 \rightarrow V_2$ such
that 
\[\gamma (Y_1(u,(x,\varphi))v) = Y_2(\gamma(u),(x,\varphi))\gamma(v)
\quad \mbox{for} \quad u,v \in V_1 ,\]
$\gamma(\mbox{\bf 1}_1) = \mbox{\bf 1}_2$, and $\gamma(\tau_1) =
\tau_2$.  

\begin{defn} {\em A} ($N = 1$ Neveu-Schwarz) vertex operator
superalgebra over $\Lambda_*$  and without odd variables {\em is a
$\frac{1}{2}${\bbb Z}-graded (by weight) $\Lambda_\infty$-module which
is also $\Bbb{Z}_2$-graded (by sign)
\[V = \coprod_{n \in \frac{1}{2} \Bbb{Z}} V_{(n)} = \coprod_{n \in
\frac{1}{2} \Bbb{Z}} V_{(n)}^0 \oplus \coprod_{n \in \frac{1}{2}
\Bbb{Z}} V_{(n)}^1 = V^0 \oplus V^1 \]  
such that only the subspace $\Lambda_*$ of $\Lambda_\infty$ acts
non-trivially on $V$, and
\[\dim V_{(n)} < \infty \quad \mbox{for} \quad n \in \frac{1}{2}
\Bbb{Z} ,\]
\[V_{(n)} = 0 \quad \mbox{for $n$ sufficiently small} ,\]
equipped with a linear map $V \otimes V \rightarrow V[[x,x^{-1}]]$, or
equivalently, 
\begin{eqnarray*}
V & \rightarrow & (\mbox{End} \; V)[[x,x^{-1}]]\\
v & \mapsto & Y(v,x) = \sum_{n \in \Bbb{Z}} v_n x^{-n-1}  
\end{eqnarray*}
where $v_n \in (\mbox{End} \; V)^{\eta(v)}$ for $v$ of homogeneous
sign in $V$, $x$ is an even formal variable, and where $Y(v,x)$ 
denotes the} vertex operator associated with $v$, {\em and equipped
also with two distinguished homogeneous vectors $\mbox{\bf 1} \in
V_{(0)}^0$ (the {\em vacuum}) and $\tau \in V_{(\frac{3}{2})}^1$.  The
following conditions are assumed for $u,v \in V$:
\[u_n v = 0 \quad \mbox{for $n \in \frac{1}{2} \Bbb{Z}$
sufficiently large;} \]
\[Y(\mbox{\bf 1}, x) = 1 \quad \mbox{(1 on the right being
the identity operator);} \]
the} creation property {\em holds:
\[Y(v,x) \mbox{\bf 1} \in V[[x]] \qquad \mbox{and} \qquad \lim_{x
\rightarrow 0} Y(v,x) \mbox{\bf 1} = v ; \]
the} Jacobi identity {\em holds: 
\pagebreak
\[x_0^{-1} \delta \left( \frac{x_1 - x_2}{x_0}\right) Y(u,x_1)Y(v,x_2)
- (-1)^{\eta(u)\eta(v)} x_0^{-1} \delta \left( \frac{ - x_2 +
x_1}{x_0} \right)Y(v,x_2)Y(u,x_1) \]
\[ \hspace{2in} = x_2^{-1} \delta \left( \frac{x_1 - x_0}{x_2} \right)
Y(Y(u,x_0)v,x_2) ;\]
the Neveu-Schwarz algebra relations hold:
\begin{eqnarray*}
\left[L(m),L(n) \right] &=& (m - n)L(m + n) + \frac{1}{12} (m^3 - m)
\delta_{m + n , 0} c , \\
\left[ G(m + \frac{1}{2}),L(n) \right] &=& (m - \frac{n - 1}{2} ) G(m
+ n + \frac{1}{2}) ,\\
\left[ G (m + \frac{1}{2} ) , G(n - \frac{1}{2} ) \right] &=& 2L(m +
n) + \frac{1}{3} (m^2 + m) \delta_{m + n , 0} c ,  
\end{eqnarray*} 
for $m,n \in \Bbb{Z}$, where 
\[G(n + \frac{1}{2}) = \tau_{n + 1} \qquad \mbox{for} \; n \in \Bbb{Z}
, i.e. \quad Y(\tau,x) = \sum_{n \in \Bbb{Z}} G (n + \frac{1}{2}) x^{-
n - \frac{1}{2} - \frac{3}{2}}, \]  
and $c \in \Bbb{C}$;
\[L(0)v = nv \quad \mbox{for} \quad n \in \frac{1}{2} \Bbb{Z} \quad
\mbox{and} \quad v \in V_{(n)} ; \]
\[ \frac{\partial}{\partial x} Y(v, x) = Y(L(-1)v,x) .\]} 
\end{defn}
The superalgebraic vertex operator superalgebra just defined is
denoted by $(V,Y(\cdot,x),\mbox{\bf 1},\tau)$.  

A consequence of the definition is that
\begin{equation}\label{G bracket}
[G(- \frac{1}{2}), Y(v,x) ] = Y(G(- \frac{1}{2})v,x) .
\end{equation}

Note that our definition of vertex operator superalgebra (without
formal variables) is an extension of the usual notion of vertex operator
superalgebra (cf. \cite{T}, \cite{DL}, and \cite{KW}) in that $V$ is a
$\Lambda_\infty$-module instead of just a vector space over {\bbb C}.

Let $(V_1, Y_1(\cdot,x),\mbox{\bf 1}_1,\tau_1)$ and $(V_2,
Y_2(\cdot,x),\mbox{\bf 1}_2,\tau_2)$ be two vertex operator 
superalgebras over $\Lambda_*$.  A {\it homomorphism} of vertex
operator superalgebras without odd formal variables is a doubly graded
$\Lambda_*$-module homomorphism $\gamma : V_1 \rightarrow V_2$ such
that 
\[\gamma (Y_1(u,x)v) = Y_2(\gamma(u),x)\gamma(v) \quad \mbox{for}
\quad u,v \in V_1 ,\] 
$\gamma(\mbox{\bf 1}_1) = \mbox{\bf 1}_2$, and $\gamma(\tau_1) =
\tau_2$.  

Let $c$ be a complex number, $\bf{SV_1}(c,*)$ be the category of 
vertex operator superalgebras over $\Lambda_*$ with odd formal
variables, and $\bf{SV_2}(c,*)$ be the category of vertex operator
superalgebras over $\Lambda_*$ without odd formal variables. 

\begin{prop}\label{superalgebras} 
For any $c \in \Bbb{C}$, the two categories $\bf{SV_1}(c,*)$ and
$\bf{SV_2}(c,*)$ are isomorphic. 
\end{prop}

{\it Sketch of proof}: We define $F_1 : \bf{SV_1}(c,*) \rightarrow
\bf{SV_2}(c,*)$ and $F_2 : \bf{SV_2}(c,*) \rightarrow \bf{SV_1}(c,*)$
by  
\begin{eqnarray*}
F_1(V,Y(\cdot,(x,\varphi)), \mbox{\bf 1}, \tau) = (V,Y(\cdot,(x,0)),
\mbox{\bf 1}, \tau), \quad \mbox{and} \quad F_1(\gamma) = \gamma \\ 
F_2(V,Y(\cdot,x), \mbox{\bf 1}, \tau) =
(V,\tilde{Y}(\cdot,(x,\varphi)), \mbox{\bf 1}, \tau), \quad \mbox{and}
\quad F_1(\gamma) = \gamma
\end{eqnarray*}
where $\tilde{Y}(v,(x,\varphi)) = Y(v,x) + \varphi Y(G(-
\frac{1}{2})v, x)$.  Using the consequences (\ref{Y with x and phi})
and (\ref{L derivative}) of the definition of $(V,
Y(\cdot,(x,\varphi)), \mbox{\bf 1}, \tau)$ and the consequence (\ref{G
bracket}) of the definition of $(V,Y(\cdot,x), \mbox{\bf 1}, \tau)$,
it is easy to see that $(V,Y(\cdot,(x,0)), \mbox{\bf 1}, \tau)$ and
$(V,\tilde{Y}(\cdot,(x,\varphi)), \mbox{\bf 1}, \tau)$ are vertex
operator superalgebras without and with odd variables, respectively.
Then clearly $F_1$ and $F_2$ are functors, and we have $F_1 \circ F_2
= 1_{\bf{SV_2}(c)}$ and $F_2 \circ F_1 = 1_{\bf{SV_1}(c)}$. $\pfbox$   

\section{The isomorphism between the category of supergeometric vertex
operator superalgebras and the category of vertex operator
superalgebras} 

The main result in \cite{H2} states that the category of geometric
vertex operator algebras and the category of vertex operator algebras
are isomorphic.  The main result that we announce is the following
analogous result in the super case.

\begin{thm}\label{supergeometric and superalgebraic} 
For any $c \in \Bbb{C}$, the two categories $\bf{SV_1}(c,*)$ and
$\bf{SG}^+(c,*)$ are isomorphic (and hence $\bf{SV_2}(c,*)$ and
$\bf{SG}^+(c,*)$ are isomorphic).   
\end{thm}

{\it Sketch of proof}: We first define a functor $F_{SG^+} :
\bf{SG}^+(c,*) \rightarrow \bf{SV_1}(c,*)$.  Given a supergeometric 
vertex operator superalgebra $(V, \nu)$ with rank $c$, define the
vacuum $\mbox{\bf 1}_\nu \in \bar{V}$
by 
\[\mbox{\bf 1}_\nu = \nu_0 (\mbox{\bf 0}) ;\]
an element $\tau_\nu \in \bar{V}$ by
\[\tau_\nu = \frac{\partial}{\partial \epsilon} \nu_0 (\mbox{\bf 0},
\{0, - \epsilon, 0,0,0... \}) ; \]
the vertex operator $Y_\nu (v_1, (x, \varphi)) = \sum_{n \in \Bbb{Z}}
(v_1)_n x^{-n-1} + \varphi \sum_{n \in \Bbb{Z}} (v_1)_{n -
\frac{1}{2}} x^{-n-1}$ associated with $v_1 \in V$ by
\[(v_1)_n v_2 + \theta (v_1)_{n - \frac{1}{2}} v_2 = \mbox{Res}_z
\left( z^n \nu_2 ((z,\theta); \mbox{\bf 0}, (1, \mbox{\bf 0}), (1,
\mbox{\bf 0}) ) (v' \otimes v_1 \otimes v_2) \right) , \]
where $\mbox{Res}_z$ means taking the residue at the singularity $z =
0$, i.e., taking the coefficient of $z^{-1}$.  

It is easy to see from the sewing and grading axioms that $\mbox{\bf
1}_\nu, \tau_\nu \in V$, and in fact, \\
$(V, Y_\nu(\cdot, (x, \varphi)), \mbox{\bf 1}_\nu, \tau_\nu)$ is a
vertex operator superalgebra with rank $c$.  The functor $F_{SG^+}$ is
defined by  
\[F_{SG^+}(V,\nu) = (V, Y_\nu (\cdot, (x, \varphi)), \mbox{\bf 1}_\nu,
\tau_\nu) \quad \mbox{and}\quad F_{SG^+}(\gamma) = \gamma .\]

We next define a functor $F_{SV}^+ : \bf{SV}_1 (c, *) \rightarrow
\bf{SG}^+(c,*)$.  Given a vertex operator superalgebra $(V, Y(\cdot,
(x,\varphi)), \mbox{\bf 1}, \tau)$ with rank $c$ and positive square
root structure, we want to define maps $(\nu_n^+)^Y : SK^+(n)
\rightarrow S{\cal F}_V (n)$, $Q \mapsto (\nu_n^+)^Y(Q)$.  For a
supersphere with three tubes $Q = ((z,\theta); (A^{(0)},M^{(0)}),
(a_0^{(1)}, A^{(1)},M^{(1)}), (a_0^{(2)}, A^{(2)},M^{(2)}))$, we
define $(\nu_2^+)^Y(Q)$ by 

\[(\nu_2^+)^Y(Q) (v' \otimes v_1 \otimes v_2) = \langle e^{- \sum_{j
\in \Z} ( A^{(0)}_j L'(j) \; +  \; M^{(0)}_{j - \frac{1}{2}} G'(j -
\frac{1}{2}) )} v', \hspace{2.2in} \] 
\[Y( e^{- \sum_{j \in \Z} (A^{(1)}_j L(j) \; + \; M^{(1)}_{j -
\frac{1}{2}} G(j - \frac{1}{2}) )} \cdot (a^{(1)}_0)^{-L(0)} \cdot
v_1, (x_1, \varphi_1)) \cdot \]
\[\hspace{2.2in} \left. e^{- \sum_{j \in \Z} (A^{(2)}_j L(j) \; + \; 
M^{(2)}_{j - \frac{1}{2}} G(j - \frac{1}{2}) )} \cdot
(a^{(2)}_0)^{-L(0)} \cdot v_2 \rangle \right|_{(x_i,\varphi_i) = 
(z_i,\theta_i)} .\]  
For elements $Q \in SK^+(n)$, $n \neq 2$, we define $(\nu_n^+)^Y(Q)$
similarly using correlation functions of products of vertex operators,
appropriately interpreting the expressions as supermeromorphic
superfunctions.  The pair $(V, (\nu^+)^Y)$ is a supergeometric vertex 
operator superalgebra with positive square root structure.  The proof
of this fact uses Propositions \ref{Gamma and Psi} and
\ref{convergence of t-series}. The functor $F_{SV}^+$ is defined by  
\[F_{SV}^+ (V, Y(\cdot, (x,\varphi)), \mbox{\bf 1}, \tau) = (V,
(\nu^+)^Y) \quad \mbox{and} \quad F_{SV}^+ (\gamma) = \gamma .\]

In fact, $F_{SG^+} \circ F_{SV}^+ = 1_{\bf{SV}_1(c,*)}$ and
$F_{SV}^+  \circ F_{SG^+} = 1_{\bf{SG}^+(c,*)}$. $\pfbox$ 

By Proposition \ref{supergeometric isomorphism} and Theorem
\ref{supergeometric and superalgebraic}, we can define the
isomorphisms of categories $F_{SG^-} = J^- \circ F_{SG^+} :
\bf{SG}^-(c,*) \rightarrow \bf{SV}_1(c,*)$ and $F_{SV}^- =
F_{SG^-}^{-1}$.  In particular,

\begin{cor}
For any $c \in \Bbb{C}$, the two categories $\bf{SV_1}(c,*)$ and
$\bf{SG}^-(c,*)$ are isomorphic (and hence $\bf{SV_2}(c,*)$ and
$\bf{SG}^-(c,*)$ are isomorphic).  Moreover, the spin structure
symmetry between $\bf{SG}^+(c,*)$ and $\bf{SG}^-(c,*)$ defines an
automorphism $J_1$ of the category $\bf{SV_1}(c,*)$ and an
automorphism $J_2$ of the category $\bf{SV_2}(c,*)$ given by  
\[J_1 (V, Y(\cdot, (x, \varphi)), \mbox{\bf 1}, \tau) = (V, Y^-(\cdot,
(x,\varphi)), \mbox{\bf 1}, - \tau) \quad \mbox{and} \quad J_1
(\gamma) = \gamma, \] 
where $Y^-(v,(x,\varphi)) = Y(v, (x,-\varphi))$, and 
\[J_2 (V, Y(\cdot, x), \mbox{\bf 1}, \tau) = (V, Y(\cdot, x),
\mbox{\bf 1}, - \tau) \quad \mbox{and} \quad J_2(\gamma) = \gamma ,\]
respectively.  In particular, we have the following commutative
diagram of categories and isomorphisms
\[\begin{array}{ccccc}
\bf{SG}^+(c,*) & \stackrel{F_{SG^+},F_{SV}^+}{\longleftrightarrow} & 
\bf{SV_1}(c,*) & \stackrel{F_1,F_2}{\longleftrightarrow} & 
\bf{SV_2}(c,*) \\
\Big\updownarrow\vcenter{%
     \rlap{$J^\pm$}} & &
\Big\updownarrow\vcenter{%
     \rlap{$J_1$}}  & & 
\Big\updownarrow\vcenter{%
     \rlap{$J_2$}}  \\ 
\bf{SG}^-(c,*) & \stackrel{F_{SG^-},F_{SV}^-}{\longleftrightarrow} & 
\bf{SV_1}(c,*) & \stackrel{F_1,F_2}{\longleftrightarrow} & \; \;
\bf{SV_2}(c,*) \; \; .
\end{array} \]
\end{cor}

Acknowledgments: I would like to thank my advisors James Lepowsky and
Yi-Zhi Huang for their support, encouragement and invaluable advice.
This work was supported in part by an American Association of
University Women Educational Foundation Fellowship.

\end{document}